\documentclass{pnastwo}

\usepackage{graphicx}
\usepackage{amsfonts,amsmath}

\usepackage{hyperref}

\copyrightyear{2012}
\issuedate{Issue Date}
\volume{Volume}
\issuenumber{Issue Number}

\newcommand{\ep}{\epsilon}

\begin{document}

\title{The averaging principle}

\author{Gadi Fibich\affil{1}{Tel Aviv University, Israel},
Arieh Gavious\affil{2}{Ono Academic College and Ben Gurion University, Israel},
Eilon Solan\affil{1}{}}

\contributor{Submitted to Proceedings of the National Academy of Sciences
of the United States of America}

\maketitle

\begin{article}
\begin{abstract}
Typically, models with a heterogeneous property are
considerably harder to analyze than the corresponding homogeneous models, in which
the heterogeneous property is replaced with its average value.
In this study
we show that any outcome of a heterogeneous model that satisfies the two
properties of \emph{differentiability} and \emph{interchangibility}, 
is~$O(\epsilon^2)$ equivalent to the outcome of the corresponding homogeneous
model, 
where~$\epsilon$ is the level of heterogeneity. We then use this  \emph{averaging principle}
to obtain new results in queueing theory, game theory (auctions), and social networks (marketing).
\end{abstract}




Mathematical modeling is a powerful tool in scientific research. Typically,
the mathematical model is merely an approximation of the actual problem.
Therefore, when choosing the model to work with, one has to strike a balance
between complex models that are more realistic, and simpler models that
are more amenable to analysis and simulations.
This dilemma arises, for example, when the model contains a heterogeneous
quantity. In such cases, a huge simplification is usually achieved by
replacing the heterogeneous quantity with its average value. The natural
question that arises is whether this approximation is ``legitimate'', i.e.,
whether the error that is introduced by this approximation is sufficiently
small.

Let us illustrate this with the following example, which will be discussed
in details later on. 
Consider a queue with $k$~heterogeneous
servers, whose expected service times are $\{ \mu_1, \dots, \mu_k\}$. We would
like to calculate analytically the expected number of customers in the
system, which we denote by~$F(\mu_1, \dots, \mu_k)$.\footnote{%
In order to focus on the heterogeneous property, we suppress the dependence
of~$F$ on other parameters.} While an explicit expression for~$F(\mu_1,
\dots, \mu_k)$ is not available, there is a well-known explicit expression
in the case of $k$~homogeneous servers, which we denote by
$F_{\mathrm{homog.}}({\mu}):= F(\underbrace{{\mu},\dots, {\mu}}_{\times k
})$.
A natural approximation for the expected number of customers in
the system is
\begin{equation}  \label{eq:approx}
F(\mu_1, \dots, \mu_k) \approx F_{\mathrm{homog.}}({\bar \mu}),
\end{equation}
where $\bar \mu$ is the average of $\{ \mu_1, \dots, \mu_k\}$.

More generally, let $F(\mu_1,\dots, \mu_k)$ denote an ``outcome'' of a
heterogeneous model, let
\begin{equation}
   \label{eq:epsilon}
\epsilon =\frac{\max_{1 \le i \le k}|\mu _{i}-\bar{\mu}|}{|\bar{\mu}|},
\end{equation}
 denote the level of heterogeneity of $\{ \mu_1, \dots, \mu_k\}$,
  and let $F_{\mathrm{homog.}}({\mu})$ denote
the outcome of the corresponding homogeneous model. If the function~$%
F(\mu_1, \dots, \mu_k)$ is differentiable, then it immediately follows that
\begin{equation*}
F(\mu_1,\dots, \mu_k) = F_{\mathrm{homog.}}(\bar{\mu}) + O(\epsilon).
\end{equation*}
Therefore, roughly speaking, for a $10\%$ heterogeneity level, the error of
approximating $F(\mu_1,\dots, \mu_k)$ with $F_{\mathrm{homog.}}(\bar{\mu%
})$ is~O(10\%). In many studies in different fields, however, researchers
have noted that the error of this approximation is considerably smaller than
$O(\varepsilon)$. Moreover, this observation seems to hold even when the
level of heterogeneity is not small.

In this study we show that these observations follow from a general
principle, which we call the \textit{Averaging Principle}. Specifically, we
show that any outcome of a heterogeneous model that satisfies the two
properties of \emph{differentiability} and \emph{interchangeability}, is~$%
O(\epsilon^{2})$ asymptotically equivalent to the outcome of the
corresponding homogeneous model, i.e.,
\begin{equation*}
F(\mu_1,\dots, \mu_k) = F_{\mathrm{homog.}}(\bar{\mu}) + O(\epsilon^2).
\end{equation*}
Thus, if the function~$F$ is also interchangeable, the error of the
approximation~\eqref{eq:approx} for a $10\%$ heterogeneity level is
only~O(1\%).


\section{The Averaging Principle}
\label{sec:principle}

Let $F(\mu_1,\dots, \mu_k)$ be an outcome of a model with a heterogeneous property,
captured by the $k$ parameters $\{\mu_1,\dots, \mu_k\}$,
that satisfies the following two properties:
\begin{enumerate}
\item \textbf{Differentiability}: $F$ is  twice-differentiable
at and near the diagonal $\mu_1 = \dots = \mu_k$.

\item \textbf{Interchangeability}: For every $(\mu_1, \dots, \mu_k) \in
\mathbb{R}^k$ and every $i \not= j$, $F(\dots, \mu_i,\dots,\mu_j,\dots)=
F(\dots, \mu_j,\dots,\mu_i,\dots)$. Thus, the outcome~$F$ is independent of the identities/indices of the heterogeneous parameters.\footnote{For example, in the queueing-system example, 
switching the identities/locations of two servers does not affect the expected number of customers in the system.}
\end{enumerate}
Then, we have the following result:\footnote{ This and all other proofs are given in the Appendix.}
\begin{theorem}[The Averaging Principle]
\label{thm:averaging}
Let $F$ satisfy the differentiability and
interchangeability properties. Let ${\mbox{\boldmath $\mu $}} =(\mu_1,
\dots, \mu_k)$ be ``sufficiently close to the diagonal'', i.e.,
\begin{equation*}
||{\mbox{\boldmath $\mu $}}-\bar{\mbox{\boldmath $\mu $}}_A||<C_{{\bar \mu}_A},
\end{equation*}
where ${\bar{\mbox{\boldmath $\mu $}}}_A = (\underbrace{\bar{\mu}_A, \dots, \bar{%
\mu}_A}_{\times k})$, $\bar{\mu}_A = \frac1k\sum_{j=1}^k \mu_j$ is the arithmetic average,
$||\cdot||$ is a vector norm on~$\mathbb{R}^k$,
and $C_{{\bar \mu}_A}$ is a positive constant that only depends on~$\bar{\mu}_A$
(and of course on $F$). Then,
\begin{equation}
\label{eq:averaging-thm}
F(\mu_1,\dots, \mu_k) = F_{\mathrm{homog.}}(\bar{\mu}_A)
+ O(||\mbox{\boldmath $\mu $}-\bar{\mbox{\boldmath $\mu $}}_A||^2),
\end{equation}
where $ F_{\mathrm{homog.}}({\mu}):= F(\underbrace{{\mu},\dots, {\mu}}_{\times k}) $.
\end{theorem}


Theorem~\ref{thm:averaging} remains valid if~$\{ \mu_1, \dots, \mu_k\}$ are
functions and not scalars, see the Game Theory (auctions) example below.

\subsection{Using the geometric and harmonic averages}

In Theorem~\ref{thm:averaging} we averaged the $\{ \mu_i\}$s using the
arithmetic mean. It is well-known in homogenization theory
that in some cases the correct homogenization is provided by the geometric or the
harmonic mean. To address the question of the ``correct''
averaging, we recall the following result:

\begin{lemma}
\label{lem:3-means} Let $\mu>0$, and let $\{h_1, \dots ,h_k\} \in \mathbb{R}$%
. Then, as $\epsilon \to 0$, the arithmetic, geometric and harmonic means of~%
$\{{\mu}+ \epsilon h_1, \dots ,{\mu}+ \epsilon h_k\}$ are~$O(\epsilon^2)$
asymptotically equivalent.
\end{lemma}

\vspace{0mm} \parindent=0pt \textbf{Proof.} \hspace{3mm} \parindent=3ex  We
can prove this result using the averaging principle. Let~$\bar{\mu}_{\mathrm{%
A}}$ denote the arithmetic mean of~$\{{\mu}+ \epsilon h_1, \dots ,{\mu}+
\epsilon h_k\}$. The geometric mean $\mu_G({\mu}+ \epsilon h_1, \dots ,{\mu}%
+ \epsilon h_k) = \left(\prod_{i=1}^k ({\mu}+ \epsilon h_i) \right)^{1/k}$
satisfies the interchangeability and differentiability properties.
Therefore, application of Theorem~\ref{thm:averaging} gives
\begin{equation*}
\mu_G({\mu}+ \epsilon h_1, \dots ,{\mu}+ \epsilon h_k) =  \mu_G(\bar{\mu}_{%
\mathrm{A}}, \dots ,\bar{\mu}_{\mathrm{A}}) + O(\epsilon^2) = \bar{\mu}_{%
\mathrm{A}} + O(\epsilon^2).
\end{equation*}
The proof for the harmonic mean~$\bar{\mu}_{\mathrm{H}} =
k/(\frac{1}{\mu_1}+\cdots+\frac{1}{\mu_k})$ is similar. $\Box$
\vspace{5mm} \parindent=3ex

From Lemma~\ref{lem:3-means} and the differentiability of~$F_{\mathrm{%
homog.}}(\mu)$ it follows that
\begin{eqnarray*}
F_{\mathrm{homog.}}(\bar{\mu}_{\mathrm{A}})  &=& F_{\mathrm{homog.}}(%
\bar{\mu}_{\mathrm{G}}) + O(\epsilon^2)
\\ &=& F_{\mathrm{homog.}}(\bar{\mu}%
_{\mathrm{H}}) + O(\epsilon^2).
\end{eqnarray*}

\begin{corollary}
Let $\mu>0$. Then, the averaging principle (Theorem~\ref{thm:averaging})
remains valid if we replace  the arithmetic mean with the geometric or
harmonic means.
\end{corollary}

A natural question is which of the three averages is ``optimal'', in the
sense that it minimizes the constant in the $O(||\mbox{\boldmath $\mu $}-%
\bar{\mbox{\boldmath $\mu $}}||^2)$~error term. The answer to this question
is model specific. It can be pursued by calculating explicitly the $O(\epsilon^2)$ term,
as we will do later on.

\subsection{Weak Interchangeability}

To extend the scope of the averaging principle, we define a
weaker interchangeability property.\footnote{See the social-networks application below for an example of a
weakly-interchangeable outcome which is not interchangeable.}
\begin{itemize}
\item[2A.] \textbf{Weak interchangeability}:
For every $\mu,\tilde \mu$ and every $1 \leq  j_0 \leq k$, if
$\mu_{j_0} = \tilde{\mu}$ and $\mu_j = \mu$ for all $j \not= j_0$,
then $F(\mu_1,\ldots,\mu_k)$ is independent of the value of~$j_0$.
\end{itemize}
Thus, $F(\mu_1,\cdots,\mu_k)$ is weakly interchangeable if,
whenever all but one of the parameters are identical,
the outcome~$F$ is independent of the identity (coordinate) of the heterogeneous parameter.

Every interchangeable function~$F$ is also weakly interchangeable,
but not vice versa.  Nevertheless, the proof of
Theorem~\ref{thm:averaging} implies that:
\begin{corollary}
The averaging principle (Theorem~\ref{thm:averaging}) remains valid if we
replace the assumption of interchangeability with the assumption of weak
interchangeability.
\end{corollary}

%

\section{Queuing theory application: An M/M/k queue with heterogeneous
service rates}

\label{sec:que}

Consider a system with $k$~servers. Server~$i$ has a random service
time that is distributed according to an exponential distribution with rate~%
$\mu _{i}$. 
Customers arrive randomly according to a Poisson distribution with arrival
rate~$\lambda$. An arriving customer is randomly allocated to one of the
non-busy servers, if such a server exists. Otherwise, the customer joins a
waiting queue, which is unbounded in length.
Once a customer is allocated to a server, he gets the service he needs and
then leaves the system. This setup is known in the Queuing
literature as M/M/k model.\footnote{%
For an introduction to queueing theory, see e.g.,~\cite{Gross-85}.} Examples
for such multi-server queuing systems are call centers, queues in banks,
parallel computing, and communications in ISDN protocols.

Let $F(\mu _{1}, \dots ,\mu _{k})$ denote the expected number of customers in
the system (i.e., waiting in the queue or receiving service) in steady
state. In the case of two heterogeneous servers, $F(\mu _{1},\mu _{2})$ can be explicitly
calculated 
(see Appendix):
\begin{lemma}
\label{lem:MMk=2} Consider an M/M/2 queue with heterogeneous servers.
The expected number of customers in the system  is given by
\begin{equation}
   \label{eq:R_nu1_nu2}
F(\mu _{1},\mu _{2})=\frac{1}{(1-\rho )^{2}}\frac{1}{\frac{1}{\rho}\frac{2
\mu_1 \mu_2}{(\mu_1+\mu_2)^2} +\frac{1}{1-\rho }} , \qquad \rho := \frac{%
\lambda}{\mu_1+\mu_2}.
\end{equation}
\end{lemma}

Finding an explicit solution for~$F(\mu _{1},\dots,\mu _{k})$ when $k \ge 3$ is computationally challenging, because
it involves solving a system of $2^{k}-1$~linear equations.
In the homogeneous case $\mu _{1}=\cdots =\mu _{k}=\mu$, however, it
is well-known that
\begin{equation}
F(\underbrace{\mu ,\dots,\mu }_{\times k})=\frac{\frac{\left(
\lambda /\mu \right) ^{k}}{k!}\frac{\frac{\lambda }{k\mu }}{1-\frac{\lambda
}{k\mu }}}{\sum_{n=0}^{k-1}\frac{\left( \lambda /\mu \right) ^{n}}{n!}+\frac{%
\left( \lambda /\mu \right) ^{k}}{k!}\frac{1}{1-\frac{\lambda }{k\mu }}}%
\frac{1}{1-\frac{\lambda }{k\mu }}+\frac{\lambda }{\mu }.  \label{eq:Rsym}
\end{equation}

The function~$F(\mu _{1}, \dots ,\mu _{k})$ can be written as a sum of solutions of a system of linear equations with
coefficients that depend smoothly on~$\mu _{1}, \dots ,\mu _{k}$ (see the appendix).
Therefore, $F$~is differentiable.
Since customers are randomly allocated to the free servers, renaming the servers does not affect the expected number of customers in
the system. Hence, $F$~is also interchangeable.
Therefore, we can use the averaging principle to obtain an explicit $O(\epsilon^{2})$ approximation for~$%
F(\mu_{1},\dots,\mu _{k})$:
\begin{theorem}
\label{lem:queueing}
Consider an M/M/k queue with heterogeneous servers whose service rates are $\{ \mu _{1}, \dots ,\mu _{k}\}$.
The expected number of customers in the system  is given by
\begin{equation*}
F(\mu _{1},\dots,\mu _{k})=F_{\mathrm{homog.}}(\bar{\mu})+O(\varepsilon ^{2}),
\end{equation*}
where $F_{\mathrm{homog.}}(\bar{\mu}):=F(\underbrace{\bar{\mu} ,\dots,\bar{\mu}}_{\times k} )$ is
given by~\eqref{eq:Rsym}, $\bar{\mu} := \frac{1}{k} \sum_{i=1}^k \mu_i$, and~%
$\epsilon $ is give by~\eqref{eq:epsilon}.
\end{theorem}
%

For example, by Theorem~\ref{lem:queueing}, the expected number of customers
with $2$~heterogeneous servers is
\begin{equation}
 \label{eq:k=2-queue-expansion}
F(\mu_1 ,\mu_2 )=F(\bar{\mu} ,\bar{\mu})+O\left(\varepsilon ^{2}\right)=\frac{4\lambda \bar{\mu} }{4\  \bar{\mu}^{2}-\lambda ^{2}}%
+O\left( \varepsilon ^{2}\right),
\end{equation}
where
\begin{equation*}
\bar{\mu} =\frac{\mu _{1}+\mu _{2}}{2},\qquad \varepsilon =\frac{\mu
_{2}-\mu _{1}}{2}.
\end{equation*}%
Indeed, substituting $\mu _{1,2}=\bar{\mu} \pm \epsilon$ in the exact expression~\eqref{eq:R_nu1_nu2} and expanding in~$%
\varepsilon$ gives~\eqref{eq:k=2-queue-expansion}.


In the case of $k=8$ heterogeneous servers, even writing the system of $%
2^{8}-1=255$ equations for the 255~unknowns is a formidable task, not to
mention solving it explicitly. By the averaging principle, however,
\begin{equation*}
F(\mu _{1},\dots,\mu _{8}) = F_{\mathrm{homog.}}(\bar{\mu})+O\left( \varepsilon ^{2}\right),
\end{equation*}
where $F_{\mathrm{homog.}}(\bar{\mu}):=F(\underbrace{\bar{\mu},\dots,\bar{\mu}}_{\times 8 })$ is given by~\eqref{eq:Rsym} with~$%
k=8$.
We ran stochastic simulations  of an M/M/8 queuing system with $8$~heterogeneous servers
using the ARENA simulation software, and used it to calculate the expected number of customers in the system. The simulation
parameters were
\begin{equation*}
\lambda = \frac{28}{\mbox{hour}}, \qquad \mu = \frac{5}{\mbox{hour}}, \quad
\mu _{i}=\mu +\varepsilon
h_{i},\quad i=1,\ldots 8,
\end{equation*}
\begin{equation*}
 (h_{1},\ldots ,h_{8}) =
(1,1.5,2,3,3.5,-2.5,-4,-4.5) \frac{1}{\mbox{hour}},
\end{equation*}%
and $\varepsilon$ varies between~$0$ and~$1$ in increments of~$0.05$.
Because $\sum_{i=1}^{k}h_{i}=0$, the average service rate is $\bar{\mu} =\mu= 5$.
Therefore, by Theorem~\ref{lem:queueing},
\begin{equation*}
F(\mu +\varepsilon h_{1}, \dots, \mu +\varepsilon
h_{8}) = F_{\mathrm{homog.}}(5)+O\left(
\varepsilon ^{2}\right).
\end{equation*}
In addition, by Equation~\eqref{eq:Rsym}, $F_{\mathrm{homog.}}(5)= 6.2314$.
%
%
%

\begin{figure}[tbp]
\begin{center}
\scalebox{.6}{\includegraphics{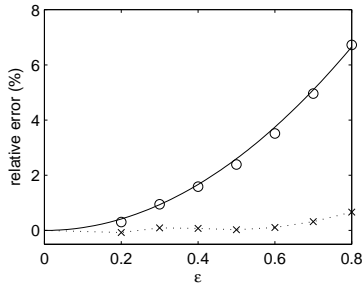}}
\end{center}
\caption{The relative error 
of the averaging-principle approximation for the steady-state number of customers
in a system with 8~heterogeneous servers, as a function of the heterogeneity
parameter~$\epsilon$. The solid line is
$\mbox{\it error} = 0.594 \epsilon^2$. The crosses denote the relative error
of the improved approximation~\eqref{eq:F_improved}.
The dotted line is $\mbox{\it error} = 0.074 \epsilon^3$.
\label{fig:queue-simulation}
}
\end{figure}

To illustrate the accuracy of this approximation,
we plot in Fig.~\ref{fig:queue-simulation} the relative error of the averaging-principle
approximation $\frac{F(\mu +\varepsilon h_{1}, \dots, \mu +\varepsilon
h_{8})-F_{\mathrm{homog.}}(5)}{F(\mu +\varepsilon h_{1}, \dots, \mu +\varepsilon
h_{8})}$.
As expected, this error scales as~$\varepsilon^{2}$.
Note that even when the
heterogeneity is not small, the averaging-principle approximation is quite
accurate. This is because the coefficient (0.594) of the $O(\epsilon^2)$~term
is small.\footnote{This coefficient will be computed analytically later on from eq.~\eqref{eq:alpha_k=8}.}
For example, when $\epsilon=0.5$ the relative error is $\approx 2\%
$, and for~$\epsilon =1$ it is below~$10\%$.

\vspace{0mm} \parindent=0pt \textbf{Remark.} \hspace{0mm} \parindent=3ex We
can also use the averaging principle to obtain $O(\epsilon^2)$%
~approximations of other quantities of interest that satisfy the interchangeability
property, such as the average waiting
time in the queue, or the probability that there are exactly $m$~customers
in the queue.

\section{Game theory application: Asymmetric Auctions}

\label{sec:auctions}

Consider a sealed-bid first-price auction with $k$~bidders,
in which the bidder who places the highest bid
wins the object and pays his bid, and all other bidders pay nothing.\footnote{%
See~\cite{Krishna} for an introduction to auction theory.}
%
A common assumption in auction theory is that of independent \emph{private-value auctions%
}, which says that each bidder knows his own valuation for the object,
does not know the valuation of the other bidders, but does know the
cumulative distribution functions (CDF) of the valuations of the other
bidders. %
Bidders are also characterized by their attitude towards risk: The
literature usually assumes that bidders are risk neutral, since this
simplifies the analysis. More often than not, however, bidders are
risk averse.

A \emph{strategy} of bidder~$i$ is a function~$b_i(\cdot)$ that assigns a bid~$b_i(v_i)$
to each possible valuation~$v_i$ of that bidder. The bid that a bidder
places depends on his valuation~$v_i$, and on his beliefs about the
distributions of the valuations of the other bidders and about their bidding
behavior. An \emph{equilibrium} in this setup is a vector of $k$~strategies $%
\{ b_i(\cdot) \}_{i=1}^k$, such that no single bidder can profit by
deviating from his bidding strategy, whatever his valuation might be, so
long that all other bidders follow their equilibrium bidding strategies.

Most of the auction literature focuses on the \emph{symmetric} (homogeneous)
case, in which the beliefs of any bidder about any other bidder (e.g., about
his distribution of valuations, his attitude towards risk, etc.) are the
same. In this case, one can look for 
%
a symmetric equilibrium, in which all bidders adopt the same strategy. In
practice, however, bidders are usually asymmetric (heterogeneous), both in
their attitude towards risk and in the distribution of their valuations. Each
bidder then faces a different competition. As a result, the equilibrium
strategies of the bidders are not the same.

The addition of asymmetry usually leads to a huge complication in the
analysis. For example, in the case of a first-price auction for a single
object with risk-neutral bidders that have private values that are
independently distributed in the unit interval~$[0,1]$ according to a common
function~$F(v)$, 
the symmetric Nash equilibrium inverse bidding strategy~$v(b) = b^{-1}(v)$
satisfies the ordinary differential equation (ODE)\footnote{%
Since we consider the case where all bidders use the same strategy, we omit
the subscript~$i$ from~$b$ and~$v$.}
\begin{equation*}
v^{\prime }(b) = \frac{1}{k-1}\frac{F(v(b))}{F^{\prime }(v(b))}\frac{1}{%
v(b)-b},\qquad v(0)=0.
\end{equation*}
This equation can be solved explicitly, yielding
\begin{equation}  \label{eq:b(v)-symm}
b(v)=v-\frac{\int_0^v F^{k-1}(s)\,ds}{F^{k-1}(v)}.
\end{equation}
Therefore, this case is ``completely understood".  From the
seller's point of view, a key property of an auction is his expected
revenue. In the symmetric case, the expression~\eqref{eq:b(v)-symm} can be
used to calculate the seller's expected revenue $R_{\mathrm{homog.}}[F]$, 
yielding
\begin{equation}  \label{eq:R_sym}
R_{\mathrm{homog.}}[F] = 1+(k-1) \int_0^1 F^k(v) \, dv - k \int_0^1 F^{k-1}(v)
\, dv.
\end{equation}

In the asymmetric case, where the value of bidder~$i$ is independently
distributed in~$[0,1]$ according to~$F_i(v)$, the inverse equilibrium
strategies $\{ v_{i}(\cdot) \}_{i=1}^k$ are the solutions of the system of ODE's
\begin{subequations}
\label{eq:intro_dvi_with_BC}
\begin{equation}  \label{eq:intro_dvi}
v_{i}^{\prime }(b)=\frac{F_{i}(v_{i}(b))}{F_{i}^\prime(v_{i}(b))}\left[
\left( \frac{1}{k-1}\sum \limits_{j=1}^{k}\frac{1}{\left( v_{j}(b)-b\right) }%
\right) -\frac{1}{\left( v_{i}(b)-b\right) }\right], 
\end{equation}
for $ i=1,\cdots,k$,
subject to the initial conditions
\begin{equation}  \label{eq:intro_dvi_left_BC}
v_i(b=0)=0,\qquad \quad i=1,\cdots,k,
\end{equation}
and the ``end condition'' at some unknown~$\bar{b}$
\begin{equation}  \label{eq:intro_dvi_right_BC}
v_i(\bar{b})=1,\quad \qquad i=1,\cdots,k.
\end{equation}
Thus, the addition of asymmetry leads to a huge complication of the
mathematical model: instead of a single ODE that can be explicitly integrated,
the mathematical model consists of a system of coupled nonlinear ODE's with a non-standard
boundary condition. As a result, the system~\eqref{eq:intro_dvi_with_BC}
cannot be explicitly solved, and it is poorly
understood, compared with the symmetric case.

In~\cite{Fibich-Gavious-2003}, Fibich and Gavious considered the system~%
\eqref{eq:intro_dvi_with_BC} in the weakly-asymmetric case $F_i = F +
\epsilon H_i$, $i=1, \dots, k$.  After several pages of
perturbation-analysis calculations, they obtained $O(\epsilon^{2})$%
~asymptotic approximations of the inverse equilibrium strategies $\{
v_{i}(b; \epsilon) \}_{i=1}^k$. Substituting these approximations in the
expression for the seller's expected revenue, showed that it is given by
\end{subequations}
\begin{eqnarray}  \label{eq:R^1st}
&&R[F_1=F + \epsilon H_1, \dots,F_k=F + \epsilon H_k] =  R_{\mathrm{homog.}}[F] \\
&& -\epsilon
(k-1)\int_{0}^{1}(1-F(v))F^{k-2}(v)\sum \limits_{i=1}^{k}H_{i}(v)\,dv
+O(\epsilon^{2}). \qquad  \notag
\end{eqnarray}

Subsequently, Lebrun~\cite{Lebrun-2009}
proved that the function on the left-hand-side of~\eqref{eq:R^1st} is differentiable in~$\ep$,
and used that to show that Eq.~\eqref{eq:R^1st} holds.
This is, in fact, a special case of the averaging principle.
Indeed, interchangeability holds since changing
the indices of the bidders  does not affect the revenue,
and, as mentioned above, differentiability in~$\epsilon$ was proved in~\cite%
{Lebrun-2009}. Therefore, by the averaging principle for functions (see the appendix),
\begin{equation}
  \label{R[F1,...,Fk]}
R[F_1=F + \epsilon H_1, \dots,F_k=F + \epsilon H_k]  = R_{\mathrm{homog.}}[\bar{F}]+O(\epsilon^{2})
\end{equation}
where $\bar{F} = F + \frac{\epsilon}{k} \sum_{i=1}^k H_i$. Substituting $\bar{F}$
in~\eqref{eq:R_sym} and expanding in powers of~$\epsilon$ gives
\begin{eqnarray*}
&& R_{\mathrm{homog.}}[\bar{F}] = R_{\mathrm{homog.}}[F]
 \\ && \quad  - \epsilon
(k-1)\int_{0}^{1}(1-F(v))F^{k-2}(v)\sum \limits_{i=1}^{k}H_{i}(v)\,dv+O(\epsilon^{2}).
\end{eqnarray*}
Hence, relation~\eqref{eq:R^1st} follows.

Numerical calculations (\cite[Table~1]{Fibich-Gavious-2003}
and~\cite[Tables~1-2]{Fibich-Gavious-Sela-94})
show that the error of the averaging-principle approximation~\eqref{R[F1,...,Fk]}
is small (typically below 1\%),
even when the asymmetry level is mild (e.g., $\epsilon=0.4$).
This provides another illustration that the averaging-principle approximation
can be useful even when $\epsilon$ is not very small.

 The averaging principle does not only lead to a simpler derivation
of relation~\eqref{eq:R^1st},  but also enables us to derive a more general novel result:
\begin{theorem}
Consider an {\em anonymous} auction\footnote{i.e., an auction in which the winner and the amount that each bidder pays depend solely on their bids,
and not on the identity of the bidders.}
in which all $k$ bidders have the same attitude towards risk,
and all bidders follow
the same ``rules'' when they determine their bidding strategies.\footnote{ For example, bidders may use
\emph{bounded rationality}~\cite{Selten-88} when determining their bidding
strategies. Thus, bidders may restrict themselves to a class of simple
strategies, such as low-order polynomial functions of the valuation~$v$.
They may even not be aware of the concept of equilibrium. Nevertheless, as
long as all bidders have the ``same'' bounded rationality, the
interchangeability requirement holds.}
Let $F_1,\cdots,F_k$ be the cumulative distribution functions of the valuations of the bidders,
and let $R[F_1,\cdots,F_k]$ be the expected revenue of the seller.
If $R$ is twice differentiable at and near the diagonal,
then
\begin{eqnarray*}
&& R[F_1,\cdots,F_k] = R_{\mathrm{homog.}}[\bar{F}] + O(\ep^2),
\end{eqnarray*}
where  $R_{\mathrm{homog.}}[\bar{F}] = R[\bar F,\cdots,\bar F]$, $\bar F$ is the average of~$F_1,\cdots,F_k$,
and $\ep$ is the level of heterogeneity.
\end{theorem}

Indeed, the assumptions of the theorem imply that~$F$ is interchangeable. Therefore, if~$F$ is also differentiable, the theorem follows from the averaging principle.

\section{Social-networks application: Diffusion of new products}

Diffusion of new products is a fundamental problem in Marketing, which has been
studied in diverse areas such as retail service, industrial technology,
agriculture, and educational, pharmaceutical and consumer-durables markets~%
\cite{Mahajan_Muller_Bass_1993}. Typically, the diffusion process begins
when the product is first introduced into the market, and progresses through a
series of adoption events. An individual can adopt the product due to
\emph{external influences} such as mass-media or commercials, and/or due to
\emph{internal influences} by other individuals who have already adopted the
product (word of mouth).
  The internal influences depend on the underlying social-network
structure, since adopters can only influence people that they ``know''.
The social network is usually modeled by an undirected graph, where each
vertex is an individual, and two vertices are connected by an edge if they
can influence each other.

The first quantitative analysis of diffusion of new products was the {\em Bass
model}~\cite{Bass_1969}, which inspired a huge body of theoretical and
empirical research. In this model and in many of the subsequent product-diffusion models:
\begin{enumerate}
\item A new product is introduced at time~$t=0$.

\item Once a consumer adopts the product, he remains an adopter at all later
times.

\item If consumer~$j$ has not adopted before time~$t$, the probability that
he adopts the product in the time interval $[t,t+s)$, given that the product
was already adopted by~$n_j(t)$ people that are connected to~$j$,
and that no other consumer adopts the
product in the time interval~$[t,t+s)$, is
\begin{eqnarray}
&& \hspace{-1cm} \text{Prob}\left( j\text{ adopts in }[t,t+s)~\Big|~n_j(t),
   \begin{array}{c} \text{ no other consumer} \\ \text{adopts in }[t,t+s)
    \end{array}
\right)   \notag \\
&&\qquad =\left( p_j+{\frac{n_j(t)}{m_j}}\cdot q_j \right) s+O(s^{2}),
\label{Bass_Formula-individual}
\end{eqnarray}%
as $s\rightarrow 0$, where $m_j$ is the total number of individuals connected to consumer~$j$ and the parameters~$p_j$ and~$q_j$ describe the likelihood of individual~$j$ to adopt the product due to
external and internal influences, respectively.
\end{enumerate}

\begin{figure}
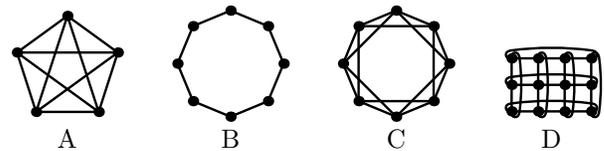

\centerline{
\includegraphics{average1.eps}
\ \ \ \ \
\includegraphics{average2.eps}
\ \ \ \ \
\includegraphics{average3.eps}
\ \ \ \ \
\includegraphics{average4.eps}
}
\caption{Examples of translation-invariant networks.  \label{fig:networks}}
\end{figure}

We say  that a  social network is \emph{translation invariant},
if any individual sees exactly the same network structure.
Therefore, in particular, $m_j$~is independent of~$j$.
Examples of translation-invariant social networks are (see Fig.~\ref{fig:networks}):
\begin{itemize}
\item[A)] A complete graph, in which any two individuals are connected.

\item[B)] A one-dimensional circle, in which each individual is connected to his two nearest neighbors.

\item[C)] A one-dimensional circle, in which each individual is connected to his four nearest neighbors.

\item[D)] A 2-dimensional torus, in which each individual is connected to his four
nearest neighbors.
\end{itemize}

We say that all individuals are homogeneous when all individuals share the same parameters, i.e.,
$p_j = p$ and $q_j=q$ for every individual~$j$.
Let~$N(t)$ denote the number of adopters at time~$t$.
The expected aggregate adoption curve~$E_{\rm homog.}[N(t;p,q)]$
in several translation-invariant social networks with homogeneous individuals were
analytically calculated in~\cite{Niu_2002,Fibich-Gibori-2010}.
In these studies, the assumption that all individuals are homogeneous was essential
for the analysis.

One of the fundamentals of marketing theory is that
consumers are anything but homogeneous.
An explicit calculation of the expected aggregate adoption curve $E[N(t;
\{p_j \},\{q_j \})]$ in the heterogeneous case, however, is much
harder than in the homogeneous case. As a result,
the effect of heterogeneity is not well understood.

The averaging principle allows us to approximate the
heterogeneous model with the corresponding homogeneous model.
Consider a translation-invariant network. Then, for~$t \geq 0$ the function $F(\{p_j\},\{q_j\}):= E[N(t; \{p_j \},\{q_j \})]$ is differentiable and weakly-interchangeable (see Appendix).
Therefore, by the averaging principle,
\begin{theorem}
   \label{lem:marketing}
The expected aggregate adoption curve in a translation invariant social network
with heterogeneous individuals, can be approximated with
$$
E[N(t;
\{p_j \},\{q_j \})] = E_{\rm homogeneous}[N(t;\bar{p},\bar{q})]+O(\epsilon^2),
$$
where 
$\bar{p}$ and $\bar{q}$ are the averages of~$\{p_j\}$ and~$%
\{q_j\}$, respectively, and~$\epsilon$ is the level of heterogeneity of~$\{p_j \}$
and~$\{q_j \}$.
\end{theorem}

Theorem~\ref{lem:marketing} is consistent
with previous numerical findings:
\begin{itemize}
  \item In~\cite{Goldenberg_Libai_Muller_2001_a}, simulations of an agent-based model
with a complete graph showed that heterogeneity in~$p$ and~$q$ had a minor effect
on the expected aggregate adoption curve.

\item Simulations of agent-based models with~1D and~2D translation-invariant networks~\cite[%
Figure 18]{Fibich-Gibori-2010} showed that when the values of~$\{p_j\}$ and~$%
\{q_j\}$ are uniformly distributed within $\pm20\%$ of the corresponding
values~$\bar{p}$ and~$\bar{q}$ of the homogeneous individuals, the
heterogeneous and homogeneous adoption curves are nearly indistinguishable.
Even when the heterogeneity level was increased to $\pm50\%$, the two
adoption curves were still very close.
\end{itemize}

\section{Calculating the $O(\protect \epsilon^2)$ term}

The averaging principle is based on a two-term Taylor expansion of~$F$.
Therefore, the error of this approximation is given, to leading order,
by the quadratic term in this expansion.
When~$F$ satisfies the differentiability and interchangeability properties\footnote{Here we cannot assume that $F$ is only weakly interchangeable,
since we require that $\frac{\partial^2 F}{\partial \mu_i \partial\mu_j}  =\frac{\partial^2 F}{\partial \mu_1 \partial\mu_2}$ for all $i,j$.}
%
and~$\mu$ is the arithmetic mean, this error is given by (see the appendix):

\begin{subequations}
   \label{eq:epsilon2}
\begin{eqnarray}
F(\mu_1, \dots, \mu_k) - F(\bar{\mu}_A, \dots, \bar{\mu}_A) \sim  \alpha
\sum_{i=1}^k (\mu_i-\bar{\mu}_A)^2,
\end{eqnarray}
where
\begin{eqnarray}
\label{eq:alpha}
 \alpha: = \frac{1}{2}\left(\frac{\partial^2 F}{\partial \mu_1 \partial\mu_1}\bigg|_{%
\bar{\mbox{\boldmath $\mu $}}_A}-\frac{\partial^2 F}{\partial \mu_1 \partial\mu_2}\bigg|_{%
\bar{\mbox{\boldmath $\mu $}}_A} \right).
\end{eqnarray}
\end{subequations}
 Therefore, 
\begin{enumerate}
  \item The magnitude of this error is~$\sim |\alpha| \,||\mbox{\boldmath $\mu $}-\bar{\mbox{\boldmath $\mu $}}_A||^2 $.
  \item The sign of this error is the same as the sign of~$ \alpha$. 

\end{enumerate}

A Taylor expansion in~$h$ gives
\begin{eqnarray*}
&& \hspace{-0.5cm} F(\bar{\mu}_A+2h,\bar{\mu}_A,\underbrace{\bar{\mu}_A, \dots,\bar{\mu}_A}_{\times k-2 })-
  F(\bar{\mu}_A+h,\bar{\mu}_A+h,\underbrace{\bar{\mu}_A \dots,\bar{\mu}_A}_{\times k-2 })
\\ && \quad \sim 2\alpha h^2, \qquad h \ll 1.
\end{eqnarray*}
This shows that in order to determine the sign of~$\alpha$, one can compare the effect
of adding $h$~units to two parameters with the corresponding effect of adding $2h$~units to a single parameter.

The value of $\alpha$ can be calculated as follows:
\begin{lemma}
   \label{lem:e2}
Assume that $F(\mu_1, \dots. \mu_k)$ satisfies the differentiability and
interchangeability properties. Then,
$$
  \alpha = \frac{k}{2(k-1)} \left( \frac{\partial^2 F}{\partial \mu_1 \mu_1}
\bigg|_{ \bar{\mbox{\boldmath $\mu $}}_A}
    -\frac{1}{k^2} F_{\mathrm{homog.}%
}^{\prime \prime }(\bar{\mu}_A) \right).
$$
\end{lemma}
Therefore, this calculation only requires the explicit calculations of~$F$ in the
homogeneous case $\mu_1 = \dots = \mu_k$,
and in the case that the heterogeneity is limited to a single coordinate
(i.e., when $\mu_2 = \dots = \mu_k$).
In many cases, this is a considerably easier task than
the explicit calculation of~$F$ in the fully-heterogeneous case.

To illustrate this, consider again the M/M/k example of Fig.~\ref{fig:queue-simulation}
with $k=8$.
While the fully-heterogeneous case requires solving
$2^{k}-1=255$ equations, the single-coordinate heterogeneous case
 requires solving only $2\cdot k =16$~equations.
Solving these 16~equations symbolically and using Lemma~\ref{lem:e2} yields, see the appendix,
\begin{equation}
  \label{eq:alpha_k=8}
\alpha(k=8)=\frac{1}{2\lambda\bar\mu} \frac{\sum_{i=0}^{12} c_i \left(\frac{\bar\mu}{\lambda}\right)^i}
{\left( \sum_{i=0}^{7} b_i \left(\frac{\bar\mu}{\lambda}\right)^i\right)^2 },
\end{equation}
where the values of $\{c_i, b_i\}$ are listed in the following table:
\begin{table}[h!]
\caption{Values of $\{c_i, b_i\}$. \label{table:bi_ci}
}
\begin{tabular}{|r||r|r|}
$i$  & $c_i$ & $b_i$  \\ \hline
0    & 1 & 1\\
1    & 45 & 14\\
2    & 999 & 126\\
3    & 14280 & 840\\
4    & 144720 & 4200\\
5    & 1088640 & 15120\\
6    & 6249600 &35280 \\
7    & 27941760 &40320 \\
8    &  97977600& \\
9    &  263390400 & \\
10    & 514382400 & \\
11   & 653184000 & \\
12    &406425600  & \\ \hline
\end{tabular}
\end{table}
In particular, substituting $\bar\mu = 5$ and~$\lambda = 28$ yields~$\alpha \approx 0.00837$. This leads to the improved approximation
\begin{eqnarray}
\nonumber
  F(\mu_1, \dots, \mu_8) &\approx& F_{\mathrm{homog.}}(\bar{\mu})+ \alpha \sum_{i=1}^{8} (\mu_{i}-\bar{\mu})^2
 \\ &\approx& F_{\mathrm{homog.}}(5)+ 0.594 \epsilon^2.
   \label{eq:F_improved}
\end{eqnarray}
The error of this improved approximation scales as $0.074\epsilon^3$,
see Fig.~\ref{fig:queue-simulation}, which is the next term in the Taylor expansion.
In particular, the relative error of~\eqref{eq:F_improved}
 is below 1.5\% for $0 \le \epsilon \le 1$.

\section{Final remarks}

The averaging principle is based on a simple observation: the leading-order effects of
 heterogeneity cancel out when the outcome is interchangeable.
 Nevertheless, it can lead to a significant simplification of mathematical models in all branches of science.
The averaging principle is unrelated to averaging that originates
from laws of large numbers in large populations, and it holds, e.g.,
when there are few servers in a queuing system, or a few bidders in
an auction.

The interchangeability and the weak interchangeability properties are
usually easy to check. The differentiability of~$F$ is easy to check in some
cases, but can be quite a challenge in others. We note, however, that more
often than not, functions that arise in mathematical models are
differentiable, unless there is a ``very good reason'' why they are not. While
this is  a very informal statement, we make it in order to point
out that the ``generic'' case is that the outcome~$F$ is differentiable, rather then the other way around.

An important issue is the ``level of heterogeneity'' that is covered by the
averaging principle. Strictly speaking, the level of heterogeneity should be
``sufficiently small''. In practice, however, in many cases the averaging
principle provides good approximations  even when $\epsilon = 0.5$.
In other words, the coefficient of the $O(\epsilon^2)$~term is~O(1).
 While
this is also an informal statement, we make it in order to point out
that one should not be ``surprised'' that the averaging principle holds even
when $\epsilon$ is not very small.


\begin{acknowledgments}
We thank Uri Yechiali for suggesting the problem of M/M/k queues with heterogeneous servers
and the solution of the M/M/2 problem.
The research of Solan was partially supported by the ISF grant \#212/09 and by the Google Inter-university center for Electronic Markets and Auctions.
\end{acknowledgments}

\newpage
\appendix

\section{Proof of Theorem~\ref{thm:averaging} }
   \label{app:proof-averagibg-thm}

Because of the differentiability of~$F$, there exists a
positive constant~$C_{\bar{\mu}_A}$, such that
for all $||\mbox{\boldmath $\mu $}- \bar{\mbox{\boldmath $\mu $}}_A||<C_{\bar{\mu}}$,
we can expand $F(\mbox{\boldmath $\mu $})$ as
\begin{equation*}
F(\mbox{\boldmath $\mu $})=F(\bar{\mbox{\boldmath $\mu $}}_A)+
\sum_{j=1}^{k}(\mu_{j}-\bar{\mu}_A)\frac{\partial F}{\partial \mu _{j}}\bigg|_{%
 \bar{\mbox{\boldmath $\mu $}}_A}+O(||\mbox{\boldmath $\mu $}- \bar{\mbox{\boldmath $\mu $}}_A||^{2}).
\end{equation*}%
Because $F$ is interchangeable,
\begin{equation*}
\frac{\partial F}{\partial \mu _{i}}\bigg|_{\bar{%
\mbox{\boldmath $\mu $}}}=\frac{\partial F}{\partial \mu _{1}}\bigg|_{\bar{\mbox{\boldmath $\mu $}}_A},\qquad j=1,\dots ,k.
\end{equation*}%
Therefore,
\begin{equation*}
F(\mbox{\boldmath $\mu $})=F(\bar{\mbox{\boldmath $\mu $}}_A)+
\frac{\partial F}{\partial \mu _{1}}\bigg|_{%
 \bar{\mbox{\boldmath $\mu $}}_A} \sum_{j=1}^{k}(\mu_{j}-\bar{\mu}_A)+O(||\mbox{\boldmath $\mu $}- \bar{\mbox{\boldmath $\mu $}}_A||^{2}).
\end{equation*}
Since $\bar{\mu}_A$ is the arithmetic average,
$\sum_{j=1}^{n}(\mu_{j}-\bar{\mu}_A)=0$. Hence, the result follows. %
%
%

\section{Proof of Lemma~\ref{lem:MMk=2} }

\begin{figure}[th]
\begin{center}
\scalebox{.7}{\includegraphics{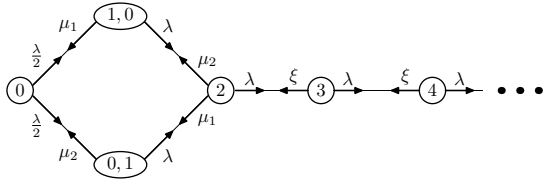}}
\end{center}
\caption{Transition diagram of a queue with two heterogeneous servers. State
\textquotedblleft 0\textquotedblright \ corresponds to the situation in which
no server is busy. States $(1,0)$ corresponds to the situation in which
server 1 is busy and server 2 is not busy. States $(0,1)$ corresponds to the
situation in which server 1 is not busy and server 2 is busy. State
\textquotedblleft $k$\textquotedblright \ for $k\geq 2$ corresponds to the
situation in which both servers are busy and $k-2$ customers wait in the
queue. \label{fig:Q2} }
\end{figure}

We calculate~$F(\mu _{1},\mu_{2})$ explicitly using the steady-state transition diagram that is shown in
Fig.~\ref{fig:Q2}.
We denote by $p_{i}$ the steady-state probability for the system to be with $%
i$~customers, and by $p_{1}^{(1,0)}$ and~$p_{1}^{(0,1)}$ the steady-state
probability for the system to be with 1~customer in server~1 and~2,
respectively. In particular, $p_{1}=p_{1}^{(1,0)}+p_{1}^{(0,1)}$. Since in
steady state the amount of inflow is equal to the amount of outflow, the
following equalities hold:
\begin{subequations}
\begin{eqnarray}
&& \lambda p_{0} = \mu _{1}p_{1}^{(1,0)}+\mu _{2}p_{1}^{(0,1)},
\label{eq:Q2_system_a} \\
&&\frac{\lambda }{2}p_{0}+\mu _{2}p_{2}=(\lambda +\mu _{1})p_{1}^{(1,0)}, \\
&&\frac{\lambda }{2}p_{0}+\mu _{1}p_{2} =(\lambda +\mu _{2})p_{1}^{(0,1)},
\label{eq:Q2_system_c} \\
&&\lambda p_{1}^{(1,0)}+\lambda p_{1}^{(0,1)}+(\mu _{1}+\mu _{2})p_{3}
=(\lambda +\mu _{1}+\mu _{2})p_{2}, \qquad\qquad   \label{eq:Q2_system_d} \\
&&\lambda p_{n}+(\mu _{1}+\mu _{2})p_{n+2} \nonumber \\
&& \qquad   =(\lambda +\mu _{1}+\mu
_{2})p_{n+1},\quad n=2,3,\dots  \label{eq:Q2_system_e}
\end{eqnarray}%
%
%
%
We can view~\eqref{eq:Q2_system_a}--\eqref{eq:Q2_system_c} as a linear
system for the three unknowns $\{p_{0},p_{1}^{(1,0)},p_{1}^{(0,1)}\}$.
Solving this system for~$p_{0}$ yields
\end{subequations}
\begin{equation*}
p_{0}=\frac{2\mu _{1}\mu _{2}}{\lambda ^{2}}p_{2}.
\end{equation*}%
In addition, the solution of~\eqref{eq:Q2_system_d}--\eqref{eq:Q2_system_e}
is $p_{n}=\left( \frac{\lambda }{\mu _{1}+\mu _{2}}\right) ^{n-2}p_{2}=\rho
^{n-2}p_{2}$ for $n\geq 1$.
Substituting the above in
\begin{equation*}
1=\sum_{n=0}^{\infty }p_{n}=p_{0}+\sum_{n=1}^{\infty }\rho
^{n-2}p_{2}=\left( \frac{2\mu _{1}\mu _{2}}{\lambda ^{2}}+\frac{1}{\rho }%
\frac{1}{1-\rho }\right) p_{2},
\end{equation*}%
gives
$p_{2}=\left( \frac{2\mu _{1}\mu _{2}}{\lambda ^{2}}+\frac{1}{\rho }\frac{1}{%
1-\rho }\right) ^{-1}.$
Therefore,
\begin{eqnarray*}
F(\mu _{1},\mu _{2}) &=&\sum_{n=0}^{\infty }np_{n}=\sum_{n=0}^{\infty }n\rho
^{n-2}p_{2}=\frac{p_{2}}{\rho }\sum_{n=0}^{\infty }n\rho ^{n-1}\\
&=&\frac{p_{2}}{%
\rho }\left( \sum_{n=0}^{\infty }\rho ^{n}\right) ^{\prime } = \frac{p_{2}}{\rho }\left( \frac{1}{1-\rho }\right) ^{\prime }=\frac{p_{2}%
}{\rho }\frac{1}{(1-\rho )^{2}},
\end{eqnarray*}%
and the result follows.

\section{M/M/3 queue}

\begin{figure}[tbp]
\begin{center}
\scalebox{.7}{\includegraphics{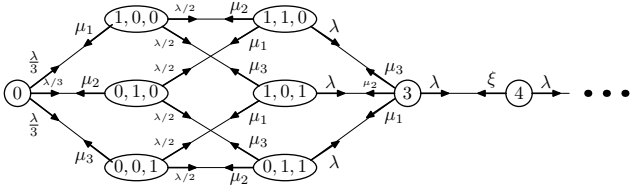}}
\end{center}
\caption{Same as Fig.~\ref{fig:Q2} with three heterogeneous servers. For
example, state $(0,1,1)$ corresponds to the situation in which server 1 is
not busy and servers~2 and~3 are busy.  \label{fig:Q3}
}
\end{figure}

Consider the case of three heterogeneous servers with average service
times $\mu _{1}$, $\mu _{2}$ and~$\mu _{3}$. Denote by $p_{0}$, $%
p_{1}^{(1,0,0)}$, $p_{1}^{(0,1,0)}$, $p_{1}^{(0,0,1)}$, $p_{2}^{(1,1,0)}$, $%
p_{2}^{(1,0,1)}$, $p_{2}^{(0,1,1)}$, $p_{3}$, $p_{4}$, \dots , the
steady-state probabilities. Thus, for example, $p_{2}^{(1,0,1)}$ is the
steady-state probability that servers~1 and~3 are busy, server~2 is free,
and there are no waiting customers in the queue (we denote by $p_n,~n \geq 2$ the
probability having $n$ customers in the system). The transition diagram for $%
k=3$~servers is given in Fig.~\ref{fig:Q3}. The steady-state
equations are
\begin{eqnarray*}
&&\lambda p_{0}=\mu _{1}p_{1}^{(1,0,0)}+\mu _{2}p_{1}^{(0,1,0)}+\mu
_{3}p_{1}^{(0,0,1)}, \\
&&\frac{\lambda }{3}p_{0}+\mu _{2}p_{2}^{(1,1,0)}+\mu
_{3}p_{2}^{(1,0,1)}=(\mu _{1}+\lambda )p_{1}^{(1,0,0)}, \\
&&\frac{\lambda }{3}p_{0}+\mu _{1}p_{2}^{(1,1,0)}+\mu
_{3}p_{2}^{(0,1,1)}=(\mu _{2}+\lambda )p_{1}^{(0,1,0)}, \\
&&\frac{\lambda }{3}p_{0}+\mu _{1}p_{2}^{(1,0,1)}+\mu
_{2}p_{2}^{(0,1,1)}=(\mu _{3}+\lambda )p_{1}^{(0,0,1)}, \\
&&\frac{\lambda }{2}p_{1}^{(1,0,0)}+\frac{\lambda }{2}p_{1}^{(0,1,0)}+\mu
_{3}p_{3}=(\lambda +\mu _{1}+\mu _{2})p_{2}^{(1,1,0)}, \\
&&\frac{\lambda }{2}p_{1}^{(1,0,0)}+\frac{\lambda }{2}p_{1}^{(0,0,1)}+\mu
_{2}p_{3}=(\lambda +\mu _{1}+\mu _{3})p_{2}^{(1,0,1)}, \\
&&\frac{\lambda }{2}p_{1}^{(0,1,0)}+\frac{\lambda }{2}p_{1}^{(0,0,1)}+\mu
_{1}p_{3}=(\lambda +\mu _{2}+\mu _{3})p_{2}^{(0,1,1)}, \\
&& \qquad =(\lambda +\mu
_{1}+\mu _{2}+\mu _{3})p_{3}, \\
&&\lambda p_{n}+\left( \mu _{1}+\mu _{2}+\mu _{3}\right) p_{n+2} \\
&& \qquad =(\lambda
+\mu _{1}+\mu _{2}+\mu _{3})p_{n+1},\qquad n\geq 3,\\
&& \sum_{n=0}^\infty p_n =1.
\end{eqnarray*}%
The solution of the last two equations is $p_{n}=\left( \frac{\lambda }{\mu
_{1}+\mu _{2}+\mu _{3}}\right)^{n-3} \!\!\!\!\!\! p_{3}$ for $n\geq 2$. The
values of $\{p_{0},p_{1},p_{2}\}$ as a function of $p_{3}$ can be evaluated
explicitly with MAPLE, by solving the first~$2^{3}-1=7$~linear equations for $\{p_{0}$, $%
p_{1}^{(1,0,0)}$, $p_{1}^{(0,1,0)}$, $p_{1}^{(0,0,1)}$, $p_{2}^{(1,1,0)}$, $%
p_{2}^{(1,0,1)}$, $p_{2}^{(0,1,1)}\}$. The resulting expression
for~$F(\mu _{1},\mu _{2},\mu _{3})$, however, is
extremely cumbersome and not informative.

\section{Proof of Theorem~\ref{lem:queueing} }

Since customers are randomly assigned to the available
servers, $F(\mu _{1},\dots ,\mu _{k})$~is interchangeable. To see that $F$
is differentiable in $(\mu _{1},\dots ,\mu _{k})$, we note that $%
F=\sum_{n=0}^{\infty }np_{n}$ where~$p_{n}$ is the steady-state probability
that there are $n$~customers in the system. In addition, $\{p_{n}\}_{n=1}^{k}
$ are the solutions of a linear system with coefficients that depend
smoothly on $(\mu _{1},\dots ,\mu _{k})$, and $p_{n}=(\frac{\lambda }{\mu
_{1}+\dots +\mu _{k}})^{n-k}p_{k}$ for $n\geq k-1$.
This was shown explicitly for the cases~$k=2$ and $k=3$; the proof for~$k>3$
is similar.

\section{Averaging principle for functions (proof of eq.~\eqref{R[F1,...,Fk]})}

Let $(F_j({\bf x})))_{j=1}^k$ be
functions in the same function space $\mathcal{F}$, and let $\epsilon \in {\Bbb R}$. Let $R : (F_1, \dots, F_k) \mapsto R[F_1, \dots, F_k] \in {\Bbb R}$ be a functional.
We say that the functional $R$ is interchangeable,
if
 $R(\dots, F_i,\dots,F_j,\dots)=
R(\dots, F_j,\dots,F_i,\dots)$ for all~$i \not=j$.
We say that the functional $R$ is differentiable if
the scalar function~$\tilde{R}(\epsilon):=R[F_1 = F+\epsilon H_1, \dots, F_k= F+\epsilon H_k]$ is twice differentiable
at and near $\epsilon=0$,
for every  $F \in \mathcal{F}$ and every $(H_j({\bf x}))_{j=1}^k \in \mathcal{F}^k$.

Given functions $(F_j({\bf x})))_{j=1}^k$ in $\mathcal{F}$, denote
$\bar{F} = \frac{1}{k}\sum_{j=1}^k F_j$ and
$H_j = F_j - \bar{F}$.
By Taylor expansion,
$$
\tilde{R}(\epsilon) = \tilde{R}(0) + \epsilon \sum_{j=1}^k \frac{\delta R}{\delta F_j} H_j
+ O(\epsilon^2),
$$
where $\frac{\delta R}{\delta F_j}$ is the variational derivative. Because $R$ is interchangeable,
$$
\tilde{R}(\epsilon) = \tilde{R}(0) + \epsilon \frac{\delta R}{\delta F_1}\sum_{j=1}^k  H_j
+ O(\epsilon^2).
$$
%
%
%
In particular, if $F = \bar{F}$, then~$\sum_{j=1}^k  H_j = 0$. Hence,
$
\tilde{R}(\epsilon) = \tilde{R}(0) + O(\epsilon^2),
$
which is~\eqref{R[F1,...,Fk]}.

\section{Proof of Theorem~\ref{lem:marketing}}

We first prove that $F$ is differentiable.
Denote $\delta_{i,i'}=1$ if individuals $i$ and $i'$ influence each other, and $\delta_{i,i'}=0$ otherwise.
For every~$k$, every set of $k$~consumers $%
\{i_{1},i_{2},\ldots ,i_{k}\}$, and every increasing sequence of times $%
0\leq t_{1}\leq \cdots \leq t_{k}$, denote by $P(i_{1},t_{1},i_{2},t_{2},%
\ldots ,i_{k},t_{k})$ the probability that consumer $i_{1}$ adopts the
product before time~$t_{1}$, consumer $i_{2}$ adopts the product between
times~$t_{1}$ and~$t_{2}$, etc., and all consumers who are not in $%
\{i_{1},\ldots ,i_{k}\}$ do not adopt the process by time~$t_{k}$. Then,
\begin{equation*}
P(i_{1},t_{1})=\big(1-\exp (-p_{i_{1}}t_{1})\big)\prod_{j\neq i_{1}}\exp
(-p_{j}t_{1}).
\end{equation*}%
Similarly,
\begin{eqnarray*}
&&P(i_{1},t_{1},i_{2},t_{2},\ldots ,i_{k},t_{k})= \\
&&\quad P(i_{1},t_{1},i_{2},t_{2},\ldots ,i_{k-1},t_{k-1}) \\
&&\times \left( 1-\exp \left( -\Big(p_{i_{k}}+\sum_{m=1}^{k-1}\delta_{i_k,i_m}q_{i_{m}}\Big)%
(t_{k}-t_{k-1})\right) \right)
 \\
&&\times \prod_{j\not\in\{i_1,\ldots,i_k\}}\exp \left( -\left(
p_{j}+\sum_{m=1}^{k-1}\delta_{j,i_m}q_{i_{m}}\right) (t_{k}-t_{k-1})\right) .
\end{eqnarray*}%
Hence, the function $P(i_{1},t_{1},i_{2},t_{2},\ldots ,i_{k},t_{k})$ is
differentiable in $\{p_{i},q_{i}\}$. Finally,
\begin{eqnarray*}
&  &E[N(t;\{p_{j}\},\{q_{j}\})] = \frac{1}{M}\sum_{\pi }\sum_{k=1}^{M}\frac{k}{(M-k!} \hfil \mbox{}
 \\
&&\hspace{-0.5cm} \times  \int_{t_{1}=0}^{t}%
\int_{t_{2}=t_{1}}^{t} \!\!\!\!\!\!\!\! \cdots
\int_{t_{k-1}=t_{k-2}}^{t} \hspace{-1cm} P(i_{1},t_{1},\ldots
,i_{k-1},t_{k-1},i_{k},t)\,dt_{k-1}\ldots dt_{1}, \qquad \mbox{}
\end{eqnarray*}%
where $\pi$ ranges over all permutations on the set of $M$ individuals.
Therefore, the differentiability of $E[N(t;\{p_{j}\},\{q_{j}\})]$
follows.

Because the network is translation invariant,
$F$ is weakly-interchangeable in $\{ p_j \}$ and in $\{q_j\}$. By this we mean that
\begin{itemize}
  \item
If $p_m = \tilde{p}$, $p_j = p$ for all $j \not= m$, and $q_j = q$ for all~$j$, then $F$ is
independent of the value of~$m$.
  \item
If $q_n = \tilde{q}$, $q_j = q$ for all $j \not= n$, and $p_j = p$ for all~$j$, then $F$ is
independent of the value of~$n$.
\end{itemize}
Therefore, the result follows from a slight modification of the proof of Theorem~\ref{thm:averaging}.

\section{Proof of equation~\eqref{eq:epsilon2}}

Since $F$ is interchangeable,
the quadratic term in the Taylor expansion of~$F(\mu_1, \dots. \mu_k)$
around  the arithmetic mean 
is equal to
\begin{eqnarray*}
&& \sum_{i,j=1}^k (\mu_i-\bar{\mu}_A) (\mu_j-\bar{\mu}_A) \frac{\partial^2 F}{\partial
\mu_i \partial \mu_j}\bigg|_{\bar{\mbox{\boldmath $\mu $}}_A} =  \\
\nonumber &&  \frac{\partial^2 F}{\partial \mu_1 \partial \mu_2}%
\bigg|_{\bar{\mbox{\boldmath $\mu $}}}
\sum_{i,j=1, i\not=j}^k (\mu_i-\bar{\mu}_A) (\mu_j-\bar{\mu}_A)
 + \frac{\partial^2 F}{\partial \mu_1 \partial \mu_1}\bigg|_{%
\bar{\mbox{\boldmath $\mu $}}}\sum_{i=1}^k (\mu_i-\bar{\mu}_A)^2.
\end{eqnarray*}
Since $\bar{\mu}_A$ is the arithmetic mean,
$$
\sum_{i,j=1}^k (\mu_i-\bar{\mu}_A) (\mu_j-\bar{\mu}_A) =
\sum_{i=1}^k (\mu_i-\bar{\mu}_A) \sum_{j=1}^k(\mu_j-\bar{\mu}_A)
= 0.
$$
Therefore, the result follows.

\section{Proof of Lemma~\ref{lem:e2}}

Consider the case where $\mu_i = \bar{\mu}+ h$ for $i=1, \dots, k$.
By equation~\eqref{eq:epsilon2},
\begin{eqnarray*}
&&\frac{1}{2} \sum_{i,j=1}^k (\mu_i-\bar{\mu}) (\mu_j-\bar{\mu}) \frac{\partial^2 F}{\partial
\mu_i \partial \mu_j}\bigg|_{\bar{\mbox{\boldmath $\mu $}}} \\
\nonumber && \quad =\frac{1}{2} \frac{\partial^2 F}{\partial \mu_1 \partial \mu_2}%
\bigg|_{\bar{\mbox{\boldmath $\mu $}}} k(k-1)h^2
 + \frac{1}{2}\frac{\partial^2 F}{\partial \mu_1 \partial \mu_1}\bigg|_{%
\bar{\mbox{\boldmath $\mu $}}} k h^2 .
\end{eqnarray*}
On the other hand, since
\begin{equation*}
F(\bar{\mu}+  h, \dots ,\bar{\mu}+  h)
= F_{\mathrm{homog.}}(\bar{\mu}+ h),
\end{equation*}
we have
\begin{equation*}
\frac{1}{2} \sum_{i,j=1}^k (\mu_i-\bar{\mu}) (\mu_j-\bar{\mu}) \frac{\partial^2 F}{\partial
\mu_i \partial \mu_j}\bigg|_{\bar{\mbox{\boldmath $\mu $}}} = \frac{h^2}{2} F_{\mathrm{homog.}%
}^{\prime \prime }(\bar{\mu}).
\end{equation*}
Therefore,
\begin{equation*}
\frac{1}{2} \frac{\partial^2 F}{\partial \mu_1 \partial \mu_2}%
\bigg|_{\bar{\mbox{\boldmath $\mu $}}} k(k-1)h^2
 + \frac{1}{2}\frac{\partial^2 F}{\partial \mu_1 \partial \mu_1}\bigg|_{%
\bar{\mbox{\boldmath $\mu $}}} k h^2  = \frac{h^2}{2} F_{\mathrm{homog.}%
}^{\prime \prime }(\bar{\mu}).
\end{equation*}
Hence,
\begin{equation*}
\frac{\partial^2 F}{\partial \mu_1 \partial \mu_2}\bigg|_{ %
\bar{\mbox{\boldmath $\mu $}}}  = \frac{1}{k-1} \left( \frac{1}{k} F_{%
\mathrm{homog.}}^{\prime \prime }(\bar{\mu}) - \frac{\partial^2 F}{%
\partial \mu_1 \partial \mu_1}\bigg|_{ \bar{\mbox{\boldmath
$\mu $}}} \right).
\end{equation*}

\section{Calculation of $\alpha$}

We illustrate the computation of the coefficient~$\alpha$ for a queue with 8~servers. Consider then the case of a single server with service time~$\mu_1$,
and seven servers with service time~$\mu$.
Denote by $p_{0,n}$ and $p_{1,n}$, $n=1, \dots, 6$, the steady-state probabilities that~$n$ out of the homogeneous
servers are busy and that the single heterogeneous servers is free or busy, respectively.
The equations for the $2\cdot 8-1 = 15$ variables $\{p_0, p_{0,1}, p_{1,0}, \dots,p_{1,6}, p_{0,7}\}$
are
\begin{eqnarray*}
&& \lambda p_{0,0}  = \mu p_{0,1}+\mu _{1}p_{1,0}\text{\ ,}  \\
&& -p_{0,0}\frac{\lambda }{8}+p_{1,0}(\lambda +\mu _{1})-p_{1,1}\mu  = 0, \\
&& p_{0,n}\left( \lambda +n\mu \right)  = p_{0,n-1}\frac{8-n}{9-n}\lambda
+p_{1,n}\mu _{1}+p_{0,n+1}(n+1)\mu,
\\
&& \quad n=1,\ldots ,6,
\label{eq:p0,n} \\
&& p_{1,n}\left( \mu _{1}+\lambda +n\mu \right) = p_{1,n-1}\lambda
+p_{1,n+1}(n+1)\mu +p_{0,n}\frac{\lambda }{8-n},
\\ && \quad n=1,\ldots ,5, \\
&& p_{0,7}(\lambda +7)\mu )=p_{0,6}\frac{%
\lambda }{2}+p_{7}\mu _{1}\rho,
\end{eqnarray*}%
where $\rho =\frac{\lambda }{7\mu +\mu _{1}}$,  $p_{n}=\rho
^{n-7}p_{7}$ for $n \ge 8$, and $\sum_{n=0}^{\infty}p_{n} = 1$.
These equations can be solved with Maple~\footnote{The Maple code is available at
\url{www.bgu.ac.il/~ariehg/averagingprinciple.html}.}, and the solution can be used to calculate
$F(\mu_1,\underbrace{\mu, \dots,\mu}_{\times 7})$ explicitly.
Differentiating this expression twice with respect to~$\mu_1$,
differentiating~$F_{\mathrm{homog.}}$, see eq.~\eqref{eq:Rsym},
twice with respect to~$\mu$,
 and using  Lemma~\ref{lem:e2}, yields equation~\eqref{eq:alpha_k=8}.
Substituting $\bar\mu=5$  and~$\lambda = 28$ gives
$\alpha \approx  0.00837$.
In addition,
 $
 \sum_{i=1}^{8} (\mu_{i}-\bar{\mu})^2  =  \epsilon^2 \sum_{i=1}^{8} h_{i}^2 = 71 \epsilon^2 .
$
Therefore,
$
\alpha \sum_{i=1}^{8} (\mu_{i}-\bar{\mu})^2 \approx 0.594 \epsilon^2.
$

\end{article}


\begin{thebibliography}{10}

\bibitem{Bass_1969} Bass, F.M., {\em A New Product Growth Model for Consumer Durables.} Management Science 15, pp. 215-227(1969).


\bibitem{Fibich-Gavious-2003} Fibich, G. and Gavious, A., {\em Asymmetric
first-price auctions - a perturbation approach,} Mathematics of
Operational Research, 28, 836–852 (2003).


\bibitem{Fibich-Gavious-Sela-94} Fibich, G. Gavious, A., and Sela, A.,
{\em Revenue Equivalence in Asymmetric Auctions,} Journal of Economic
Theory, 115, 309-321 (2004).

\bibitem{Fibich-Gibori-2010} Fibich, G. and Gibori R., {\em Aggregate
diffusion dynamics in agent-based models with a spatial structure,} Operations Research, 58,1450-1468 (2010).


\bibitem{Goldenberg_Libai_Muller_2001_a} Goldenberg, J., Libai, B. and
Muller, E., {\em Using Complex Systems Analysis to Advance Marketing Theory
Development.} Academy of Marketing Science Review. [Online]
01, special issue on Emergent and Co-evolutionary Processes in Marketing (2001).

\bibitem{Gross-85} Gross, D. and Harris, C. M., {\em A Fundamentals of
Queueing Theory}, 2nd ed., John Wiley \& Sons, New York, (1985).

%

\bibitem{Krishna} Krishna, V. {\em Auction Theory}, Academic Press, San
Diego (2002).

\bibitem{Lebrun-2009} Lebrun, B., {\em Auctions with Almost Homogeneous
Bidders,} Journal of Economic Theory, 144, 1341-1351 (2009).

\bibitem{Mahajan_Muller_Bass_1993} Mahajan, V., Muller, E. and Bass, F. M.,
{\em New-Product Diffusion Models}. In J. Eliashberg and G. L. Lilien,
editors, Handbooks in Operations Research and Management Science,
Volume 5: Marketing, pp. 349-408. North-Holland, Amsterdam, The Netherlands (1993).




\bibitem{Niu_2002} Niu, S.C., {\em A Stochastic Formulation of the Bass
Model of New-Product Diffusion}. Mathematical Problems in
Engineering 8: 249-263 (2002).

%
\bibitem{Selten-88} Selten, R., {\em Features of experimentally observed
bounded rationality,} European Economic Review, 42, 413-436 (1988).


\end{thebibliography}
\end{document}